\documentclass[11pt]{article}
\usepackage{latexsym}
\usepackage{amssymb}
\font\gorditas = msbm8
\def\bbb#1{\hbox {{\gordas #1}}}
\def\errita{\hbox{\gorditas R}}
                                                                                
\font\gordas = msbm10 at 12pt
\def\bbb#1{\hbox {{\gordas #1}}}
\def\erre{{\bbb R}}

\def\UNO{1\mkern-7mu1}

\newtheorem{theorem}{Theorem}[section]
\newtheorem{lemma}[theorem]{Lemma}

\newtheorem{proposition}[theorem]{Proposition}
\newtheorem{definition}[theorem]{Definition}

\newtheorem{remarks}[theorem]{Remarks}

\begin{document}
\noindent
\vglue.3cm

\begin{center}
{\large\bf OCCUPATION TIME FLUCTUATIONS OF AN INFINITE VARIANCE 
BRANCHING SYSTEM IN LARGE DIMENSIONS}\\[.5cm]
(OCCUPATION TIME FLUCTUATIONS IN LARGE DIMENSIONS)\\[1cm]
\end{center}
\noindent
TOMASZ BOJDECKI\\
Institute of Mathematics\\
University of Warsaw\\
ul. Banacha 2\\
02-097 Warsaw, Poland\\
E.mail: tobojd@mimuw.edu.pl\\[.5cm]
LUIS G. GOROSTIZA $^*$\\
Department of Mathematics\\
Centro de Investigaci\'on y de Estudios Avanzados\\
A.P. 14-740\\
M\'exico, 07000 D.F., Mexico\\
E.mai: lgorosti@math.cinvestav.mx\\[.5cm]
ANNA TALARCZYK\\
Institute of Mathematics\\
University of Warsaw\\
ul. Banacha 2\\
02-097 Warsaw, Poland\\
E.mail: annatal@mimuw.edu.pl

\vglue1cm
\noindent
We prove limit theorems for rescaled occupation time fluctuations of a $(d,\alpha,\beta)$-branching particle system (particles moving in $\erre^d$ according to a  spherically symmetric $\alpha$-stable L\'evy process, $(1+\beta)$-branching, $0<\beta<1$,  uniform Poisson initial state), in the cases of critical dimension, $d=\alpha(1+\beta)/\beta$, and large dimensions, $d>\alpha(1+\beta)/\beta$. The fluctuation processes are continuous but their limits are stable processes with independent increments, which have jumps. The convergence is in the sense of finite-dimensional distributions, and also of space-time random fields (tightness does not hold in the usual Skorohod topology). The results are in sharp contrast with those for intermediate dimensions, $\alpha/\beta<d<d(1+\beta)/\beta$, where the limit process is continuous and has long range dependence (this case is studied by Bojdecki et al, 2005). The limit process is measure-valued for the critical dimension, and ${\cal S}'(\erre^d)$-valued for large dimensions. We also raise some questions of interpretation of the different types of dimension-dependent results obtained in the present and previous papers in terms of properties of the particle system.
\vglue.5cm
\noindent
\newpage
\section{Introduction}
\setcounter{section}{1}
\setcounter{equation}{0}

We consider the $(d,\alpha,\beta)$-branching particle system, which consists of particles evolving independently in $\erre^d$ according to a spherically symmetric $\alpha$-stable L\'evy process (called standard $\alpha$-stable process henceforth) and a $(1+\beta)$-branching law, $0<\beta<1$. This law has generating function
\begin{equation}
\label{eq:1.1}
s+\frac{(1-s)^{1+\beta}}{1+\beta},\quad 0<s<1,
\end{equation}
is critical, and belongs to the domain of attraction of a stable law with exponent $1+\beta$ (the case $\beta=1$ corresponds to binary branching). The particle lifetime distribution is exponential with parameter $V$ (this parameter is not
 particularly  relevant in this paper, but we keep it for consistency with our  previous papers, Bojdecki et al 2004, 2006a, 2006b, 2005, where sometimes it plays a role). We assume that the system starts off at time $0$ from a Poisson random field with intensity measure $\lambda$ (Lebesgue measure).

The $(d,\alpha,\beta)$-branching particle system and its associated superprocess have been widely studied (see Dawson, 1993, Dawson et al, 1989, Dawson and Perkins, 1991, Etheridge, 2000, Fleischmann and G\"artner, 1986, Gorostiza and Wakolbinger, 1991, Iscoe, 1986, M\'el\'eard and Roelly, 1990, for some of the early results). Although this spatial branching system is  special, the different types of behaviors it exhibits depending on relationships between $d,\alpha$ and $\beta$ (e.g. the persistence/extinction dichotomy, Gorostiza and Wakolbinger, 1991) make it an interesting test case.

Let $N=(N_t)_{t\geq 0}$ denote the empirical measure process of the system, i.e., $N_t(A)$ is the number of particles in the set $A\subset \erre^d$ at time $t$. The rescaled occupation time fluctuation process is defined by
\begin{equation}
\label{eq:1.2}
X_T(t)=\frac{1}{F_T}\int^{Tt}_0(N_s-\lambda)ds=\frac{T}{F_T}\int^t_0(N_{Ts}-\lambda)ds,\quad t\geq 0,
\end{equation}
where $F_T$ is a suitable norming for convergence and $T$ is the scaling parameter which accelerates time and tends to $\infty$. Note that $EN_s=\lambda$ for all $s>0$, due to the initial Poisson condition, the criticality of the branching and the $\alpha$-stable motion.

The process $X_T$ exhibits different asymptotic behaviors as $T\rightarrow \infty$, depending on relationships between the parameters $d,\alpha,\beta$. For {\it intermediate dimensions}, $\alpha/\beta<d<\alpha(1+\beta)/\beta$, and $F_T=T^{(2-\beta-\frac{d}{\alpha}\beta)/(1+\beta)}$, the limit process has the form $K\lambda\xi$, where $K$ is a constant and $\xi=(\xi_t)_{t\geq 0}$ is a continuous, self-similar, stable process which has long range dependence (Bojdecki et al, 2005). In the present paper we consider the cases of {\it critical dimension}, $d=\alpha(1+\beta)/\beta$, and {\it large dimensions}, $d>\alpha(1+\beta)/\beta$. We prove limit theorems for $X_T$ described briefly as follows. For the critical dimension and $F_T=(T\log T)^{1/(1+\beta)}$, the limit has the form $K\lambda\xi$, where $K$ is a constant and $\xi=(\xi_t)_{t\geq 0}$ is a real standard $(1+\beta)$-stable L\'evy process totally skewed to the right. For large dimensions and $F_T=T^{1/(1+\beta)}$, the limit is an ${\cal S}'(\erre^d)$-valued stable process with stationary independent increments $({\cal S}'(\erre^d)$ is the space of tempered distributions, dual of the space ${\cal S}(\erre^d)$ of smooth rapidly decreasing functions). These two limit processes have jumps; recall that a continuous process with independent increments is necessarily Gaussian (Kallenberg, 2002, Theorem 13.4); this is in fact the case for $\beta=1$ (Bojdecki et al, 2006b).

Clearly, since the fluctuation processes are continuous and the limit processes have jumps, there cannot be convergence in the  Skorohod space $D([0,1], {\cal S}'(\erre^d))$ with the usual $J_1$-topology. We prove convergence of finite-dimensional distributions and also in the sense of space-time random fields 
(Bojdecki et al, 1986). We do not know if tightness holds in a weaker topology, e.g., the $M_1$-topology (Skorohod, 1956), or the $S$-topology 
(Jakubowski, 1997).

The main observations regarding the results for critical and large dimensions are the striking fact that the fluctuation processes $X_T$, which are continuous, develop jumps in the limit as  $T\rightarrow\infty$, and that the limit processes 
have independent increments. These properties are in sharp contrast with those for intermediate dimensions (Bojdecki et al, 2005), where, moreover, convergence takes place in the space of continuous functions $C([0,1],{\cal S}'(\erre^d)$).

In Bojdecki et al (2006a, 2006b) we proved functional limit theorems for the system with $\beta=1$, which corresponds to finite variance (binary) branching, for intermediate, critical and large dimensions. In that case all the limits are continuous. Another relevant paper is Birkner and Z\"ahle (2005), where functional limit theorems are proved for fluctuations of the occupation time of the origin for a system of critical binary branching random walks on the d-dimensional lattice,  $d\geq 3$, with methods different from ours; the results for  initial Poisson state are parallel to those of Bojdecki et al (2006a, 2006b) and 
 and in addition the equilibrium case is treated. In contrast with Birkner and Z\"ahle (2005), we have also investigated the spatial structure in our case.
 For the $(d,\alpha,\beta)$-branching particle system with $d=\alpha/\beta$, there is a functional ergodic theorem (Talarczyk, 2005). The present paper has partial origins in Iscoe (1986), where the occupation time of the $(d,\alpha,\beta)$-superprocess is studied. The occupation times of the $(d,\alpha,\beta)$-branching particle system and the $(d,\alpha,\beta)$-superprocess have analogous properties, but the particle system is technically more involved.

The methods of proof in this paper are similar to those in 
Bojdecki et al (2006a, 2006b, 2005), with some different technical complexities; the Fourier transform tools used before for $\beta=1$ are not applicable with $\beta <1$. 
We refer to Bojdecki et al (2005) for some technical points.
 We remark that
the critical dimension is more difficult to deal with than the large 
dimensions, although the limit process is simpler.

Additive processes in ${\cal S}'(\erre^d)$ and other nuclear spaces are an interesting subject on its own right (e.g. 
It\^o, 1980, P\'erez-Abreu et al, 2005, \"Ustunel, 1984). This article shows that processes of that type actually arise in physical models.

Section 2 contains the results, and Section 3 the proofs. In Section 4 we make some comments and  raise questions of interpretation of the results for the intermediate, critical and large dimensions in terms of the particle system; these are questions that require further research.

\section{Results}
\setcounter{section}{2}
\setcounter{equation}{0}

The following notation will be used. $\langle\quad,\quad\rangle$ denotes pairing of spaces in duality (e.g., ${\cal S}'(\erre^k)$ and ${\cal S}(\erre^k))$. $||\quad||_p$ stands for $L^p$-norm, and it will be clear from the context whether the underlying space is $\erre^d$ or an interval of $\erre$. Constants are written $C,C_1,\ldots,$ with possible dependencies in parenthesis. $\Rightarrow$ denotes convergence in law in appropriate spaces, and convergence of processes in the sense of finite-dimensional distributions is denoted 
$\mathop{\Longrightarrow}\limits_{f}$.

We also need another, less known, notion of convergence, related to space-time random fields, introduced in Bojdecki et al (1986). To any stochastic process $X=(X(t))_{t\geq 0}$ with paths in the Skorohod space of cadlag functions $D(\erre_+,{\cal S}'(\erre^d))$ and any $\tau>0$, there corresponds an ${\cal S}'(\erre^{d+1})$-valued random element $\widetilde{X}$ defined by 
\begin{equation}
\label{eq:2.1}
\langle\widetilde{X},\Phi\rangle=\int^\tau_0\langle X(t),\Phi(\cdot, t)\rangle dt,\quad\Phi\in{\cal S}(\erre^{d+1}).
\end{equation}

\begin{definition}
\label{d:2.1}
{\rm
Let $X,X_n,n=1,2,\ldots,$ be ${\cal S}'(\erre^d)$-valued cadlag processes. We say that the laws of $X_n$ converge to $X$ in the space-time, or integral, sense (denoted $X_n\mathop{\Longrightarrow}\limits_{i} X$) if for each $\tau>0$,
$$\widetilde{X}_n\Rightarrow\widetilde{X}\quad{\rm as}\quad n\rightarrow\infty.$$

$\mathop{\Longrightarrow}\limits_{i}$ convergence, resembles 
$\mathop{\Longrightarrow}\limits_{f}$ convergence but none of them implies the other. Each one of these convergences together with tightness on 
$D([0,\tau],{\cal S}'(\erre^d)$) implies convergence in law in 
$D([0,\tau],{\cal S}'(\erre^d))$. 
Conversely, weak functional convergence implies 
$\mathop{\Longrightarrow}\limits_{i}$}.
\end{definition}

Let $p_t$ denote the transition density of the standard $\alpha$-stable process on $\erre^d$, and ${\cal T}_t$ the corresponding semigroup, i.e., ${\cal T}_t\varphi=p_t *\varphi$. Denote the potential operator
\begin{equation}
\label{eq:2.2}
G\varphi(x)=\int^\infty_0{\cal T}_t\varphi(x)dt=C_{\alpha,d}\int_{\errita^d}\frac{\varphi(y)}{|x-y|^{d-\alpha}}dy,
\end{equation}
where $C_{\alpha,d}=\Gamma(\frac{d-\alpha}{2})(2^\alpha\pi^{\frac{d}{2}}\Gamma(\frac{\alpha}{2}))^{-1}$ (we always assume $d>\alpha$).

The main result of this paper is the following theorem for the occupation time fluctuation process of the $(d,\alpha,\beta)$-branching particle system defined by (1.2).

\begin{theorem}
\label{t:2.2}
Assume $0<\beta<1$.

\noindent
(a) Let
\begin{equation}
\label{eq:2.3}
d>\frac{\alpha(1+\beta)}{\beta}
\end{equation}
and $F_T=T^{\frac{1}{1+\beta}}$. Then
$$X_T\mathop{\Longrightarrow}\limits_{i} X\quad{and}\quad X_T
\mathop{\Longrightarrow}\limits_{f} X\quad{as}\quad T\rightarrow\infty,$$
where $X$ is an ${\cal S}'(\erre^d)$-valued $(1+\beta)$-stable process with stationary independent increments whose distribution is determined by
\begin{equation}
\label{eq:2.4}
E{\rm exp}\{i\langle X(t),\varphi\rangle\}={\rm exp}\biggl\{-K^{1+\beta}t\int_{\errita^d}|G\varphi(x)|^{1+\beta}\biggl(
1-i({\rm sgn}G\varphi(x))\tan\frac{\pi}{2}(1+\beta)\biggr)dx\biggr\},
\end{equation}
$$\phantom{rrrrrrrrrrrrrrrrrrrrrrrrrrrrrrrrrrrrrrrrrrrrrr}
\varphi\in{\cal S}(\erre^d), t\geq 0,$$ 
where
$$K=\biggl(-\frac{V}{1+\beta}\cos\frac{\pi}{2}(1+\beta)\biggr)^{\frac{1}{1+\beta}}.$$
(b) Let
\begin{equation}
\label{eq:2.5}
d=\frac{\alpha(1+\beta)}{\beta}
\end{equation}
and $F_T=(T\log T)^{\frac{1}{1+\beta}}$. Then
$$X_T\mathop{\Longrightarrow}\limits_{i}K_1\lambda\xi\quad{and}\quad X_T
\mathop{\Longrightarrow}\limits_{f}K_1\lambda \xi\quad{as}\quad T\rightarrow\infty,$$%
where $\xi$ is a real $(1+\beta)$-stable process with stationary independent increments whose distribution is determined by
\begin{equation}
\label{eq:2.6}
E{\rm exp}\{iz\xi_t\}={\rm exp}
\biggl\{-t|z|^{1+\beta}
\left(1-i({\rm sgn}z)\tan\frac{\pi}{2}(1+\beta)\right)\biggr\},\,\,
z\in\erre,t\geq 0,
\end{equation}
and
$$K_1=\left(-V\cos\frac{\pi}{2}(1+\beta)\int_{\errita^d}
\left(\int^1_0p_r(x)dr\right)^\beta p_1(x)dx\right)^{\frac{1}{1+\beta}}.$$
\end{theorem}

\begin{remarks}
\label{r:2.3}

{\rm 
(1) Existence of a family of finite-dimensional distributions of the process $X$ described in Theorem 2.2 (a) follows, e.g., from the proof. Then the existence of the ${\cal S}'(\erre^d)$-process $X$ itself is a consequence of classical properties of nuclear spaces (regularization theorem. It\^o, 1980).

\noindent
(2) The process $\xi$ in Theorem 2.2 (b) is the standard $(1+\beta)$-stable L\'evy process totally skewed to the right.

\noindent
(3) The weak limit of $\langle X_T(t),\varphi\rangle$ in Theorem 2.2 (a) coincides with the result obtained by Iscoe (1986, Theorem 5.6) for the $(d,\alpha,\beta)$-superprocess.

\noindent
(4) For $\beta=1$ and large dimensions, in Theorem 2.2 (a) of 
Bojdecki et al (2006b) there is an additional term in the limit process which comes from the free motion of the particles (i.e., without branching). The limit is an ${\cal S}'(\erre^d)$-Wiener process, which is the sum of two independent ${\cal S}'(\erre^d)$-Wiener processes. The additional term does not appear with $\beta<1$ because in this case the branching produces larger fluctuations (of order $T^{1/(1+\beta)})$ than those of the free motion (of order $T^{1/2})$.}
\end{remarks}

It is clear that weak functional convergence of $X_T$ does not hold (with $J_1$-topology). Nevertheless, it turns out that the family $\{X_T\}_T$ posesses a property which is ``not far'' from tightness in $C([0,\tau],{\cal S}'(\erre^d))$. Namely, we have the following proposition.

\begin{proposition}
\label{p:2.4}
Assume $d\geq \frac{\alpha(1+\beta)}{\beta}$ 
and let $F_T$ be as in Theorem 2.2. Then
for any $\varphi\in {\cal S}(\erre^d)$ and $\tau>0$, 
\begin{equation}
\label{eq:2.7}
P(|\langle X_T(t_2),\varphi\rangle-\langle X_T(t_1),\varphi\rangle|>\delta)\leq \frac{C(\varphi)}{\delta}|t_2-t_1|
\end{equation}
for all $t_1,t_2\in [0,\tau],$ all $T\geq 2$ and all $0<\delta<1$.
\end{proposition}

\section{Proofs}
\label{sec:3}
\setcounter{equation}{0}

Before the proof of Theorem \ref{t:2.2} we state a simple lemma regarding the operator $G$ given by (2.2),  which will be used frequently (see, e.g., 
Iscoe, 1986, Lemma 5.3).
\begin{lemma}
\label{L3.1} If $\varphi$ is a measurable function on $\erre^d$ such that 
$$
\sup_{x\in \errita^d}(1+|x|^p) |\varphi (x)| <\infty
$$
for some $p>d$, then
$$
\sup_{x\in \errita^d} (1+|x|^{d-\alpha})|G\varphi (x)|<\infty .
$$
\end{lemma}
{\bf Proof of Theorem 2.2.} Without loss of generality we assume $\tau =1$. 

We need the form of the Laplace transform 
$E{\rm exp}\{-\langle\widetilde{X}_T,\Phi\rangle\}$ for 
$\Phi \in {\cal S}(\erre^{d+1})$, $\Phi \geq 0$, 
where $\widetilde{X}_T$ is defined by (\ref{eq:2.1}).

Denote
\begin{equation}
\label{eq:3.1}
\Psi (x,t) =\int^1_t \Phi (x,r) dr,\,\,\,\, \Psi_T (x,t)=\frac{1}{F_T}\,\,\Psi 
\biggl(x, \frac{t}{T}\biggr),
\end{equation}
and define
\begin{equation}
\label{eq:3.2}
v_T(x,t)=1-E \exp\biggl\{-\int^t_0 \langle N^x_r, \Psi_T (\cdot, T-t-r)dr\biggr\},
\end{equation}
where $N^x_r$ is the the empirical measure of the branching system started from a single particle at $x$. Identically as in 
Bojdecki et al (2006a) (see also Bojdecki et al, 2006b, Lemma 3.1), using the Feynman-Kac formula and the form of the generating function of the branching law given by (\ref{eq:1.1}), it can be shown that $v_T$ satisfies
\begin{equation}
\label{eq:3.3}
v_T(x,t)=\int^t_0 {\cal T}_{t-r} \biggl[\Psi_T (\cdot, T-r)(1-v_T (\cdot, r))-\frac{V}{1+\beta} v^{1+\beta}_T (\cdot, r)\biggr](x)dr, \quad 0\leq t\leq T,
\end{equation}
and
\begin{eqnarray}
\label{eq:3.4}
\lefteqn{
E \exp\{-\langle \widetilde{X}_T, \Phi \rangle\} =E\exp \biggl\{\int^T_0 \langle N_r, \Psi_T (\cdot, r) \rangle dr-\int_{\errita^d}\int^T_0 \Psi_T (x, r) drdx\biggr\}}\\
\label{eq:3.5}
&=&\exp \biggl\{\int_{\errita^d} \int^T_0 \Psi_T (x, T-r)v_T(x,r) drdx+\frac{V}{1+\beta}\int_{\errita^d}\int^T_0 v_T^{1+\beta}(x,r)drdx\biggr\}.
\end{eqnarray}

We will use frequently the following estimates for $v_T$:
\begin{equation}
\label{eq:3.6}
0\leq v_T (x,t)\leq 1,
\end{equation}
by (\ref{eq:3.2}), and
\begin{equation}
\label{eq:3.7}
v_T (x,t) \leq \int^t_0 {\cal T}_{t-r}\Psi_T (\cdot, T-r)(x)dr,
\end{equation}
since $1-e^{-x}\leq x$, $x\geq 0$, and 
$E\langle N^x_t, \varphi \rangle ={\cal T}_t \varphi (x)$. ((\ref{eq:3.7}) also follows from (\ref{eq:3.3}) and (\ref{eq:3.6})).

For the convergence of finite-dimensional distributions we need also the corresponding Laplace transform. 
For $\varphi_1, \varphi_2,\ldots, \varphi_k \in {\cal S}(\erre^d)$, all
$\varphi_j \geq 0$, and $0\leq t_1\leq t_2 \leq \cdots \leq t_k \leq 1$, it is easy to see that $E\exp \{-\sum^k_{j=1}\langle X_T (t_j), \varphi_j\rangle \}$ has the form (\ref{eq:3.4}) with
\begin{equation}
\label{eq:3.8}
\Psi (x,t)=\sum^k_{j=1}\varphi_j(x) \UNO_{[0,t_j]}(t).
\end{equation}
Moreover, approximating $\Psi$ by smooth functions in (\ref{eq:3.3}), we obtain that for this Laplace transform (\ref{eq:3.5}) also holds with a corresponding $v_T$ given by (\ref{eq:3.2}).

It will be convenient to write the right hand side of (\ref{eq:3.5}) in the form
\begin{equation}
\label{eq:3.9}
\exp\biggl\{I_1(T)+\frac{V}{1+\beta}(I_2 (T)-I_3 (T))\biggr\},
\end{equation}
where
\begin{eqnarray}
\label{eq:3.10}
I_1(T)&=&\int_{\errita^d} \int^T_0 \Psi_T (x, T-r) v_T (x,r) drdx,\\
\label{eq:3.11}
I_2 (T) &=&\int_{\errita^d} \int^T_0 \biggl(\int^r_0 {\cal T}_{r-u} \Psi_T (\cdot, T-u) (x) du\biggr)^{1+\beta} drdx,\\
\label{eq:3.12}
I_3 (T) &=& \int_{\errita^d} \int^T_0 \biggl[\biggl(\int^r_0 {\cal T}_{r-u} \Psi_T (\cdot, T-u)(x)du \biggl)^{1+\beta}-v^{1+\beta}_T(x,r)\biggr]drdx.
\end{eqnarray}
\vglue .5cm

{ We now prove convergence in  case (a)}.

Firstly, we want to show that
\begin{equation}
\label{eq:3.13}
\lim_{T\rightarrow \infty} E\exp{-\langle \widetilde{X}_T, \Phi \rangle} 
=\exp \biggl\{\frac{V}{1+\beta}\int_{\errita^d} \int^1_0 
(G\Psi (\cdot, r)(x))^{1+\beta}drdx\biggr\}
\end{equation}
for $\Phi \in {\cal S}(\erre^{d+1})$, $\Phi \geq 0$, and $\Psi$ given by (\ref{eq:3.1}). To simplify notation we consider $\Phi$ of the form
\begin{equation}
\label{eq:3.14}
\Phi (x,t)=\varphi (x) \psi(t),\quad \varphi \in {\cal S}(\erre^d),\,\,\, \psi \in {\cal S}(\erre), \quad\varphi, \psi \geq 0.
\end{equation}
Denote
\begin{equation}
\label{eq:3.15}
\chi (t) =\int^1_t \psi (r)dr,\quad \chi_T (t)=\chi \biggl(\frac{t}{T}\biggr),
\quad \varphi_T (x) =\frac{1}{F_T} \varphi (x).
\end{equation}

Using (\ref{eq:3.9})-(\ref{eq:3.12}), we will show that 
\begin{eqnarray}
\label{eq:3.16}
I_1 (T) &\rightarrow& 0,\\
\label{eq:3.17}
I_2 (T) &\rightarrow & \int_{\errita^d} \int^1_0 (G\varphi (x) \chi (r))^{1+\beta} drdx,\\
\label{eq:3.18}
I_3 (T) &\rightarrow & 0,
\end{eqnarray}
as $T\rightarrow \infty$.

By (\ref{eq:3.7}), (\ref{eq:3.15}) and boundedness of $\chi$ we have
\begin{eqnarray*}
I_1 (T) &\leq & C \frac{1}{F^2_T}\int_{\errita^d}\int^T_0 \varphi (x) \int^r_0 {\cal T}_u \varphi (x) dudrdx\\
&\leq & C\frac{T}{T^{\frac{2}{1+\beta}}}\int_{\errita^d}\varphi (x) G\varphi (x) dx,
\end{eqnarray*}
by (\ref{eq:2.2}). $G\varphi$ is bounded by Lemma 3.1, hence (\ref{eq:3.16}) follows since $\beta <1$.

Next, we use (\ref{eq:3.15}) and make obvious substitutions to obtain
\begin{equation}
\label{eq:3.19}
I_2 (T) =\int_{\errita^d}\int^1_0 \biggl(\int^{T(1-r)}_0 
{\cal T}_u \varphi (x) \chi \biggl(r+\frac{u}{T}\biggr)du\biggr)^{1+\beta}drdx.
\end{equation}
It is clear that 
$$
\lim_{T\rightarrow \infty} \int^{T(1-r)}_0 {\cal T}_u \varphi (x) 
\chi \biggl(r+\frac{u}{T}\biggr)du=G\varphi (x) \chi (r),
$$
hence (\ref{eq:3.17}) follows by the dominated convergence theorem and 
Lemma 3.1, since 
\begin{equation}
\label{eq:3.20}
(1+\beta) (d-\alpha)>d
\end{equation}
by (\ref{eq:2.3}).

To prove (\ref{eq:3.18}) we apply the obvious inequality
$$
a^{1+\beta}-b^{1+\beta}\leq (1+\beta)a^\beta (a-b)\quad {\rm for}\quad a\geq b\geq 0,
$$
and we obtain, by (\ref{eq:3.12}) and (\ref{eq:3.7}),
\begin{eqnarray*}
I_3 (T) &\leq & (1+\beta)\int_{\errita^d}\int^T_0 
\biggl(\int^r_0 {\cal T}_{r-u} \varphi_T (x) \chi_T (T-u) du\biggr)^\beta\\
&&\quad \times \biggl(\int^r_0{\cal T}_{r-u}\varphi_T (x) \chi_T (T-u) du-v_T (x,r)\biggr) drdx\\
&=&(1+\beta)\int_{\errita^d}\int^T_0 \biggl(\int^r_0 {\cal T}_{r-u} \varphi_T (x) \chi_T (T-u) du\biggr)^\beta\\
&&\quad \times \biggl(\int^r_0 {\cal T}_{r-u} \biggl(\varphi_T (\cdot )\chi_T (T-u) v_T (\cdot, u) 
+\frac{V}{1+\beta} v_T ^{1+\beta}(\cdot, u)\biggr)(x) du\biggr) drdx,
\end{eqnarray*}
by (\ref{eq:3.3}). We use (\ref{eq:3.7}) and boundedness of $\chi$ to arrive at
$$
I_3 (T) \leq C(J_1(T)+J_2(T)),
$$
where
\begin{eqnarray*}
J_1(T) &=& \int_{\errita^d} \int^T_0 \biggl(
\int^r_0 {\cal T}_{r-u}\varphi_T (x) du\biggr)^\beta \biggl(
\int^r_0 {\cal T}_{r-u} \biggl(\varphi_T(\cdot) \int^u_0 {\cal T}_{u-u'} \varphi_T(\cdot)du'\biggr)(x)du\biggr)drdx,\\
J_2(T) &=& \int_{\errita^d}\int^T_0 \biggl(\int^r_0 {\cal T}_{r-u} \varphi (x) du\biggr)^\beta\biggl(\int^r_0 {\cal T}_{r-u} \biggl(\int^u_0 {\cal T}_{u-u'} \varphi_T(\cdot)du'\biggr)^{1+\beta}(x)du\biggr)drdx.
\end{eqnarray*}
By (2.9) and boundedness of $G\varphi$ (Lemma 3.1) we have
$$
J_1(T)\leq CT^{-\frac{2+\beta}{1+\beta}}\int_{\errita^d}\int^T_0 (G\varphi (x))^{1+\beta}drdx\leq C_1 T^{-\frac{1}{1+\beta}},
$$
by (\ref{eq:3.20}) and Lemma \ref{L3.1}.

Similarly,
$$
J_2(T) \leq T^{-\frac{\beta}{1+\beta}}\int_{\errita^d}
(G\varphi (x))^\beta G((G\varphi)^{1+\beta})(x)dx.
$$
Hence $J_2 (T) \rightarrow 0$, since the integral above is finite by (\ref{eq:3.20}) and Lemma \ref{L3.1} applied three times (the function $(G\varphi)^{1+\beta}$ satisfies the assumption of the lemma). This finishes the proof of 
(\ref{eq:3.18}), and  therefore (\ref{eq:3.13}) is proved.

It is not difficult to see that for any $\Phi \in {\cal S}(\erre^{d+1})$, the real random variable $\langle \widetilde{X}, \Phi \rangle$ (where $X$ is the process in Theorem 2.2 (a)) has a $(1+\beta)$-stable law with characteristic function
\begin{eqnarray}
\label{eq:3.21}
\lefteqn{E\exp\{i z\langle \widetilde{X}, \Phi\rangle\}}\nonumber\\
&=&\exp \biggl\{-\,K^{1+\beta}|z|^{1+\beta}\int_{\errita^d}
\,\int^1_0 |G\Psi (\cdot, s)(x)|^{1+\beta}\left(1-i ({\rm sgn}\;z)({\rm sgn}
(G\Psi (\cdot, s)(x)))\tan \frac{\pi}{2}(1+\beta)\right)dsdx\biggr\}\nonumber\\
\nonumber\\
\end{eqnarray}
(recall that $\Psi$ is defined by (\ref{eq:3.1})).

Now, (\ref{eq:3.13}) implies that
\begin{equation}
\label{eq:3.22}
\langle \widetilde{X}_T, \Phi \rangle \Longrightarrow \langle \widetilde{X}, \Phi \rangle \quad {\rm as}\quad T\rightarrow \infty,
\end{equation}
for any $\Phi \in {\cal S}(\erre^{d+1})$. This is almost immediate for $\Phi \geq 0$, and for general $\Phi$ one employs weak convergence of two-dimensional random variables. See 
Bojdecki et al (2006a, Lemma 3.4 and Corollary 3.5) for more details (see also 
Iscoe, 1986).

By nuclearity of ${\cal S}(\erre^{d+1})$, (\ref{eq:3.22}) implies $\widetilde{X}_T \Rightarrow \widetilde{X}$ as $T\rightarrow \infty$,  hence $X_T 
\mathop{\Longrightarrow}\limits_{i} X$.

Finally, convergence of finite-dimensional distributions of $X_T$ is obtained in exactly the same way by virtue of the previous remarks, using $\Psi$ of the form (\ref{eq:3.8}).

\vglue .5cm
{Next we prove convergence in case (b).}

Again, we prove first convergence $\mathop{\Longrightarrow}\limits_{i}$. To this end we will show that
\begin{equation}
\label{eq:3.23}
\lim_{T\rightarrow \infty}E\exp\{-\langle \widetilde{X}_T,\Phi\rangle\}=
\exp\biggl\{K_2 \int^1_0 \biggl(\int_{\errita^d}\Psi (x, r)dx\biggr)^{1+\beta}dr\biggr\},
\end{equation}
where $\Psi$ is defined by (\ref{eq:3.1}), and 
$$
K_2 =V\int_{\errita^d} \biggl(\int^1_0 p_u (x) du\biggr)^\beta p_1 (x) dx.
$$

Analogously as before, we consider $\Phi =\varphi \otimes \psi$, and using (\ref{eq:3.9})-(\ref{eq:3.12}) we prove (\ref{eq:3.16}), (\ref{eq:3.18}) and
\begin{equation}
\label{eq:3.24}
I_2 (T) \rightarrow \frac{1+\beta}{V} K_2 \biggl(\int_{\errita^d}\varphi (x) dx\biggr)^{1+\beta}\int^1_0 \chi^{1+\beta}(r) dr.
\end{equation}

In spite of the simpler form of the limit, this case requires a more involved argument.

Let us start with $I_2 (T)$. Substituiting $r'=\frac{r}{T}$, then $u'=Tr'-u$, and using (\ref{eq:3.15}) we write
\begin{equation}
\label{eq:3.25}
I_2 (T) =R_1(T)+R_2 (T),
\end{equation}
where
\begin{equation}
\label{eq:3.26}
R_1(T)=\frac{1}{\log T}\int_{\errita^d}\int^1_0 
\biggl(\int^{Tr}_0 {\cal T}_u \varphi (x) \chi (1-r) du\biggr)^{1+\beta}drdx,
\end{equation}

\begin{equation}
\label{eq:3.27}
R_2(T)=\frac{1}{\log\;T}\int_{\errita^d}\int^1_0 \biggl[\biggl(\int^{Tr}_0 {\cal T}_u \varphi (x) \chi \biggl(1-r+\frac{u}{T}\biggr)du \biggr)^{1+\beta} -\biggl(\int^{Tr}_0 {\cal T}_u \varphi (x) \chi (1-r)du\biggr)^{1+\beta}\biggr]drdx.
\end{equation}

To prove (\ref{eq:3.24}) we show that $R_1 (T)$ converges to the desired limit and $R_2 (T)\rightarrow 0$.

Note that, by (\ref{eq:2.5}),
$$
\int_{\errita^d} (G\varphi (x))^{1+\beta} dx=\infty
$$
if $\varphi \not\equiv 0$, hence we can use the L'H\^opital rule to obtain 
\begin{equation}
\label{eq:3.28}
\lim_{T\rightarrow \infty} R_1 (T) =\lim_{T\rightarrow \infty}T(1+\beta) \int_{\errita^d}\int^1_0 \chi^{1+\beta}(1-r) \biggl(\int^{Tr}_0 {\cal T}_u \varphi (x) du\biggr)^\beta {\cal T}_{Tr} \varphi (x) r drdx.
\end{equation}

After subsitituting $u'=\frac{u}{Tr}$ and writing ${\cal T}_u$ in terms of $p_u$, the expression under the limit on the right hand side of (\ref{eq:3.28}) has the form
\begin{equation}
\label{eq:3.29}
(1+\beta)\int^1_0 \int_{\errita^d}\chi^{1+\beta}(1-r)\biggl(\int_0^1
\int_{\errita^d} Tr p_{uTr} (x-y) \varphi (y) dydu\biggr)^\beta 
\int_{\errita^d} Trp_{Tr}(x-z) \varphi (z) dz dx dr.
\end{equation}
We apply the self-similarity property of the stable density,
\begin{equation}
\label{eq:3.30}
p_{au}(x) =a^{-\frac{d}{\alpha}}p_u (x a^{-\frac{1}{\alpha}}),
\end{equation}
substitute $x'=x(Tr)^{-\frac{1}{\alpha}}$ and use (\ref{eq:2.5}) (which implies $(1-\frac{d}{\alpha})(1+\beta)=-\frac{d}{\alpha})$. (\ref{eq:3.29}) 
now becomes
\begin{eqnarray}
\label{eq:3.31}
\lefteqn{\kern -2cm (1+\beta)\int^1_0 \int_{\errita^d}\chi^{1+\beta}(1-r)\biggl(\int^1_0 \int_{\errita^d} p_u (x-y(Tr) ^{-\frac{1}{\alpha}})\varphi(y)dydu\biggr)^\beta \int_{\errita^d}p_1 (x-z(Tr)^{-\frac{1}{\alpha}})\varphi (z) dz dxdr}\nonumber\\
&=& (1+\beta)\int^1_0 \chi^{1+\beta}(1-r) \int_{\errita^d} 
\biggl(\biggl(\int^1_0 p_u du\biggr) * 
\widetilde{\varphi}_{Tr}(x)\biggr)^\beta p_1 * \widetilde{\varphi}_{Tr}(x) dxdr,
\end{eqnarray}
where
\begin{equation}
\label{eq:3.32}
\widetilde{\varphi}_t(y)=t^{\frac{d}{\alpha}}\varphi(yt^{\frac{1}{\alpha}}).
\end{equation}

As $\int^1_0p_udu\in L^1(\erre^d)$ we have
$$
\biggl(\int^1_0p_udu\biggr) * \widetilde{\varphi}_{Tr}\rightarrow \int^1_0p_udu\int_{\errita^d}\varphi(x)dx\quad{\rm in}\quad L^1(\erre^d),$$
and the $L^1$-norms are bounded in $r\in (0,1]$.

Analogously
$$p_1 * \widetilde{\varphi}_{Tr}\rightarrow p_1\int_{\errita^d}\varphi(x)dx\quad{\rm in}\quad L^{\frac{1}{1-\beta}}(\erre^d),$$
and the $L^{\frac{1}{1-\beta}}$-norms are bounded. By (\ref{eq:3.28}) this proves that $R_1(T)$ converges to the expression in (\ref{eq:3.24}).

Note that we have also shown that
\begin{equation}
\label{eq:3.33}
\lim_{T\rightarrow\infty}\frac{1}{\log T}\int_{\errita^d}\biggl(\int^T_0{\cal T}_u\varphi(x)du\biggr)^{1+\beta}dx=\frac{1+\beta}{V}K_2\biggl(\int_{\errita^d}\varphi(x)dx\biggr)^{1+\beta},
\end{equation}
which will be used later.

We now turn to $R_2(T)$. Note that,
by (\ref{eq:3.15}), the difference inside $[\ldots]$ in (\ref{eq:3.27}) is negative; hence, using the elementary inequality
$$(a+b)^{1+\beta}-a^{1+\beta}\leq b^{1+\beta}+(1+\beta)a^{\frac{1+\beta}{2}}b^{\frac{1+\beta}{2}},\,\,\,\, a,b\geq 0,$$
(applying it to $-[\ldots]$) we have
\begin{eqnarray*}
\lefteqn{
|R_2(T)|\leq\frac{1}{\log T}\int_{\errita^d}\int^1_0\biggl(\int^{Tr}_0{\cal T}_u\varphi(x)\biggl[\chi(1-r)-\chi\biggl(
1-r+\frac{u}{T}\biggr)\biggr]du\biggr)^{1+\beta}drdx}\\
&&\qquad +\frac{1+\beta}{\log T}\int_{\errita^d}\int^1_0\biggl(\int^{Tr}_0{\cal T}_u\varphi(x)\biggl[\chi(1-r)-\chi\biggl(1-r+\frac{u}{T}\biggr)\biggr]du\biggr)^{\frac{1+\beta}{2}}\\
&&\qquad \times\biggl(\int^{Tr}_0{\cal T}_u\varphi(x)\chi\biggl(1-r+
\frac{u}{T}\biggr)du\biggr)^{\frac{1+\beta}{2}}drdx.
\end{eqnarray*}
By the Schwarz inequality applied to the second summand we obtain
\begin{equation}
\label{eq:3.34}
|R_2(T)|\leq W(T)+(1+\beta)\sqrt{W(T)}\sqrt{I_2(T)},
\end{equation}
where $I_2(T)$ is defined by (\ref{eq:3.11}) and
\begin{equation}
\label{eq:3.35}
W(T)=\frac{1}{\log T}\int_{\errita^d}\int^1_0\biggl(\int^{Tr}_0{\cal T}_u\varphi(x)\biggl[\chi(1-r)-\chi\biggl(1-r+\frac{u}{T}\biggr)\biggr]du\biggr)^{1+\beta}drdx.
\end{equation}

Note that
$$I_2(T)\leq C\frac{1}{\log T}\int_{\errita^d}\biggl(\int^T_0{\cal T}_u
\varphi(x)du\biggr)^{1+\beta}dx,$$
so, by (\ref{eq:3.33}), $I_2(T)$ is bounded.

By (\ref{eq:3.34}), to prove that $R_2(T)\rightarrow 0$ and thus complete the proof of (\ref{eq:3.24}), it remains to show that $W(T)\rightarrow 0$.

We use the fact that $\chi$ is a Lipschitz function and we substitute $u'=\frac{u}{T}$; then
\begin{equation}
\label{eq:3.36}
W(T)\leq\frac{C}{\log T}\int_{\errita^d}\int^1_0\biggl(\int^r_0\int_{\errita^d}T p_{Tu}(x-y)\varphi(y)udydu\biggr)^{1+\beta}drdx.
\end{equation}
Applying the self-similarity property (\ref{eq:3.30}), substituting $x'=xT^{-\frac{1}{\alpha}}, y'=y T^{-\frac{1}{\alpha}}$, using (\ref{eq:2.5}) and estimating 
$\int^r_0\ldots du$ by $\int^1_0\ldots du$, we obtain
\begin{eqnarray}
W(T)&\leq&\frac{C}{\log T}\int_{\errita^d}\biggl(\int_{\errita^d}\int^1_0 p_u(x-y)uT^{\frac{d}{\alpha}}\varphi(yT^{\frac{1}{\alpha}})dydu\biggr)^{1+\beta}dx\nonumber\\
\label{eq:3.37}
&=&\frac{C}{\log T}||h *\widetilde{\varphi}_T||^{1+\beta}_{1+\beta}\leq\frac{C}{\log T}||h||^{1+\beta}_{1+\beta}||\widetilde{\varphi}_T||^{1+\beta}_1,
\end{eqnarray}
where $h(x)=\int^1_0p_u(x)udu$ and $\widetilde{\varphi}_T$ is defined by 
(\ref{eq:3.32}).

We have $||\widetilde{\varphi}_T||_1=\int_{\errita^d}\varphi(x)dx$, and it is not difficult to show that $h\in L^{1+\beta}(\erre^d)$ using self-similarity of $p_u$, the well-known estimate
$$p_1(x)\leq \frac{C}{1+|x|^{d+\alpha}},$$
and (\ref{eq:2.5}) once again. By (\ref{eq:3.37}) we obtain $W(T)\rightarrow 0$ as $T\rightarrow \infty$, and (\ref{eq:3.24}) is proved.

(\ref{eq:3.16}) follows by the same estimates as in Bojdecki et al (2005,
 see 
(\ref{eq:3.22}), (\ref{eq:3.23}) therein).

Finally, (\ref{eq:3.18}) can be obtained in exactly the same way as in 
Bojdecki et al (2005, see (\ref{eq:3.16}) and (\ref{eq:3.27})-(\ref{eq:3.36}) therein). The only difference is that under (\ref{eq:2.5}), $f(x)=\int^1_0p_u(x)du$  belongs to $L^p(\erre^d)$ for any $p<1+\beta$ (and not to $L^{1+\beta}(\erre^d)$), so the Young inequality should be applied appropriately, and then the corresponding norms of 
$g_{1,T}, g_{2,T}$ (see  Bojdecki et al 2005, (3.31)) will have the forms $C(\varphi)T^\varepsilon$ for $\varepsilon$ arbitrarily small. This finishes the proof of (\ref{eq:3.23}).

Convergence $X_T\mathop{\Longrightarrow}\limits_i Y=K\lambda\xi$ now follows as in the previous step. One should only observe that
\begin{eqnarray*}
\lefteqn{E\exp\{iz\langle\widetilde{Y},\Phi\rangle\}}\\
&=&{\rm exp}
\biggl\{
-K^{1+\beta}_1|z|^{1+\beta}\int^1_0
\biggl|
\int_{\errita^d}
\Psi(x,s)dx
\biggr|^{1+\beta}
\biggl(1-i\biggl({\rm sgn}
\biggl(z
\int_{\errita^d}\Psi(x,s)dx
\biggr)\biggr)\tan\frac{\pi}{2}(1+\beta)
\biggr)ds
\biggr\}.
\end{eqnarray*}

Also, convergence of finite-dimensional distributions can be derived similarly as before. The only difference is that now (\ref{eq:3.36}) does not hold because $\chi$ is not Lipschitz. We consider $\chi(r)=\UNO_{[0,t]}(r)$ for any fixed $t\in [0,1]$ (see (\ref{eq:3.8})). So, to prove that $W(T)$ defined by 
(\ref{eq:3.35}) tends to $0$ we argue as follows. With  this form of $\chi$ we have
$$W(T)=\frac{1}{\log T}\int_{\errita^d}\int^t_0\biggl(\int^{1-r}_{t-r}\int_{\errita^d}T p_{Tu}(x-y)\varphi(y)dydu\biggr)^{1+\beta}drdx.$$
By the self-similarity of $p_u$, substituting $x'=xT^{-\frac{1}{\alpha}},y'=yT^{-\frac{1}{\alpha}}$, and using (\ref{eq:2.5}) we obtain
$$W(T)=\frac{1}{\log T}\int_{\errita^d}\int^t_0\biggl(
\int_{\errita^d}\int^{1-t+r}_{r} p_{u}(x-y)\widetilde{\varphi}_T(y)dudy
\biggr)^{1+\beta}drdx,$$
where $\widetilde{\varphi}_T$ is given by (\ref{eq:3.32}). Hence, by the Young inequality,
$$W(T)\leq \frac{1}{\log T}\int^t_0
\biggl|\biggl|\int^{1-t+r}_r p_udu\biggr|\biggr|^{1+\beta}_{1+\beta}||\varphi||^{1+\beta}_1dr$$
(remember that $||\widetilde{\varphi}_T||_1=||\varphi||_1)$. It suffices to observe that
\begin{eqnarray*}
\lefteqn{
\int^t_0\biggl|\biggl|\int^{1-t+r}_r p_udu\biggr|\biggr|^{1+\beta}_{1+\beta}dr=\int^t_0r^{-\frac{1}{1+\beta}}\int_{\errita^d}\biggl(\int^{1-t+r}_rr p_u(x)du\bigg)^{1+\beta}dxdr}\\
&\leq&\int^t_0r^{-\frac{1}{1+\beta}}dr\int_{\errita^d}\biggl(\int^1_0 up_u(x)du\biggr)^{1+\beta}dx<\infty,
\end{eqnarray*}
by the argument following (\ref{eq:3.37}). This proves that $W(T)\rightarrow 0.$ So the convergence $X_T\mathop{\Longrightarrow}\limits_f K_1\lambda\xi$ is established.

\hfill$\Box$
\vglue.5cm
\noindent
{\bf Proof of Proposition 2.4.}
We consider $\tau=1, 0<\delta<1,0\leq t_1<t_2\leq 1$.

Arguing as in the proof of Proposition 3.3 of Bojdecki et al (2005), 
it suffices to show that
\begin{equation}
\label{eq:3.38}
P(|\langle\widetilde{X}_T,\varphi\otimes\psi\rangle|\geq\delta)\leq
\frac{C(\varphi)}{\delta}(t_2-t_1)
\end{equation}
for any $\psi\in{\cal S}(\erre)$ such that
$$\chi(t)=\int^1_t\psi(s)ds$$
satisfies
\begin{equation}
\label{eq:3.39}
0\leq\chi\leq \UNO_{[t_1,t_2]},
\end{equation}
and each $\varphi\in{\cal S}(\erre^d), \varphi\geq 0$.

Repeating the argument of that proposition (see Bojdecki et al, 2006c,
 (3.39)-(3.51)
) it is enough to prove that
\begin{equation}
\label{eq:3.40}
I\leq C(\varphi)(t_2-t_1)
\end{equation}
and
\begin{equation}
\label{eq:3.41}
II\leq C(\varphi)(t_2-t_1),
\end{equation}
where
\begin{eqnarray}
I&=&\int_{\errita^d}\int^T_0\varphi_T(x)\chi_T(T-s)\int^s_0{\cal T}_{s-r}\varphi_T(x)\chi_T(T-r)drdsdx,\\
II&=&\int_{\errita^d}\int^T_0\biggl(\int^s_0{\cal T}_{s-r}\varphi_T(x)\chi_T(T-r)dr\biggr)^{1+\beta}dsdx,
\end{eqnarray}
and $\varphi_T,\chi_T$ are given by (\ref{eq:3.15}).

We have, after obvious substitutions and by (\ref{eq:3.39}),
\begin{eqnarray*}
I&\leq&\frac{T^2}{F^2_T}
\int_{\errita^d}\int^1_0\int^s_0\varphi(x)
\chi(1-s){\cal T}_{rT}\varphi(x)drdsdx\\
&=&\frac{T^2}{F^2_T}\frac{1}{(2\pi)^d}
\int_{\errita^d}\int^1_0\int^s_0|\widehat{\varphi}(x)|^2
e^{-rT|x|^\alpha}dr\chi(1-s)dsdx,
\end{eqnarray*}
where we have used the Plancherel formula and the well known fact that 
$\widehat{{\cal T}_u\varphi}(x)=e^{-u|x|^\alpha}\widehat{\varphi}(x)$ \\
($\widehat{\,\,}$ denotes Fourier transform). Hence
\begin{eqnarray*}
I&\leq&C\frac{T}{F^2_T}\int_{\errita^d}|\widehat{\varphi}(x)|^2\frac{1}{|x|^\alpha}dx\int^1_0\chi(1-s)ds\\
&\leq&C(\varphi)\frac{T}{F^2_T}(t_2-t_1).
\end{eqnarray*}
This implies (\ref{eq:3.40}) for both $d=\frac{\alpha(1+\beta)}{\beta}$ and $d>\frac{\alpha(1+\beta)}{\beta}$.

To prove (\ref{eq:3.41}) we write
\begin{eqnarray*}
II&=&\frac{T}{F^{1+\beta}_T}
\int_{\errita^d}\int^1_0\bigg(\int^s_0T{\cal T}_{T(s-r)}\varphi(x)\chi(1-r)dr\bigg)^{1+\beta}dsdx\\
&\leq&\frac{T}{F^{1+\beta}_T}\int_{\errita^d}||f_{x,T} * g||^{1+\beta}_{1+\beta}dx,
\end{eqnarray*}
where
$$f_{x,T}(r)=T{\cal T}_{Tr}\varphi(x)\UNO_{[0,1]}(r),\quad g(r)=\chi(1-r)\UNO_{[0,1]}(r).$$

By the Young inequality we obtain
$$II\leq\frac{T}{F^{1+\beta}_T}\int_{\errita^d}||f_{x,T}||^{1+\beta}_1dx||g||^{1+\beta}_{1+\beta}.$$
Since $||g||^{1+\beta}_{1+\beta}\leq (t_2-t_1)$ by (\ref{eq:3.39}), it suffices to show that
$$\sup_{T\geq 2}\frac{T}{F^{1+\beta}_T}\int_{\errita^d}||f_{x,T}||^{1+\beta}_1dx<\infty.$$
This fact follows from Lemma 3.1 in the case $d>\frac{\alpha(1+\beta)}{\beta}$ and from (\ref{eq:3.33}) in the case $d=\frac{\alpha(1+\beta)}{\beta}$.

\hfill $\Box$

\section{Comments and questions of interpretation}
\label{sec:4}
\setcounter{equation}{0}

In the present and previous papers we have proved limit theorems for occupation time fluctuations of $(d,\alpha,\beta)$-branching particle system and described some properties of the limit processes. Some of the results raise questions concerning their meaning in terms of properties of the particle system. We mention here some of these questions which require further research.

\subsection{Transitions between intermediate and large dimensions}
\label{sub:4.1}

The results in Bojdecki et al (2005) and the present ones show that in the passage from intermediate to large dimensions the spatial structure of the limit process goes from simple $(\lambda)$, corresponding to perfect correlation in the case $\beta=1$, to complicated (truly ${\cal S}'(\erre^d)$-valued), and the temporal structure goes from complicated (long range dependence) to simple (independent increments); at the critical dimension both are simple. In subsection 4.4 we comment on the temporal change. We do not have an explanation for the spatial change. On the other hand, the size of the fluctuations of the occupation time process, measured by $F_T$, which is typically larger than for the classical central limit theorem in the case of long range dependence, does not pass continuously from intermediate to large dimensions (as $\alpha$ or $\beta$ vary); at the critical dimension they are larger by a logarithmic factor. This type of phenomenon is known to occur in some stochastic spatial models.
\subsection{Continuity and jumps}
\label{sub:4.2}
The main question about the results in this paper is to understand why the $(1+\beta)$-branching with $\beta <1$ causes the occupation time fluctuation process to generate jumps in the limit as $T\rightarrow \infty$ for the critical and large dimensions, but not  the intermediate ones.
\subsection{Poisson vs. equilibrium}
\label{sub:4.3}
The $(d, \alpha, \beta)$-branching particle system has equilibrium states in dimensions $d>\alpha /\beta$ (Gorostiza and Wakolbinger, 1991), and one may consider the system started off from  an equilibrium state instead of Poisson $(\lambda)$. Birkner and Z\"ahle (2005) give functional limit theorems with equilibrium and Poisson initial states for the fluctuations of the occupation time of the origin  of branching random walks on the $d$-dimensional lattice, $d\geq 3$. In 
Bojdecki et al (2004) we did covariance calculations with intermediate dimensions for the $(d, \alpha, 1)$-system in equilibrium. On the basis of the results of Birkner and Z\"ahle (2005), and Bojdecki et al (2004, 2006a, 2006b), it is natural to expect that for the $(d, \alpha, 1)$-system in equilibrium it is possible to  prove functional limit theorems for the occupation time fluctuations with the same normings $F_T$ as for the initial Poisson condition, and that the results would be as follows. For intermediate dimensions, $\alpha < d<2\alpha$, the limit would be of the form $K\lambda \xi$, where $\xi$ is fractional Brownian motion. For critical and large dimensions, $d\geq 2 \alpha$, the results would be the same as the corresponding ones with initial Poisson 
(Bojdecki et al, 2006b). Carrying this  further, we think that limit theorems can also be proved for the $(d, \alpha, \beta)$-system with $\beta <1$ in equilibrium. The normings would be the same as for initial Poisson, and the limits would have the same forms as in the Poisson case for critical and large dimensions. For intermediate dimensions the limit is expected to be  of the form $K\lambda \eta$, where $\eta$ is a continuous, self-similar, $(1+\beta)$-stable process with stationary increments and long range dependence, which should  be a kind of ``fractional stable process''.

Going back to $\beta=1$ and intermediate dimensions (covariance calculations in Bojdecki et al, 2004), in the equilibrium case we get, for the temporal part, fractional Brownian motion (fBm) $\xi$ with covariance function
$$
\frac{1}{2} (s^h +t^h-|s-t|^h),
$$
where $h=3-d/\alpha$, and in the Poisson case we get sub-fractional Brownian motion (sub-fBm) $\zeta$ with covariance function
$$
s^h+t^h-\frac{1}{2} [(s+t)^h+|s-t|^h],
$$
where $h=3-d/\alpha$ again. $\xi$ is defined also for $t<0$, and the two processes are related by
\begin{equation}
\label{eq:4.1}
\zeta\; \displaystyle{\stackrel{d}{=}}\; \biggl(\frac{1}{\sqrt{2}}(\xi_t+\xi_{-t})\biggr)_{t\geq 0}.
\end{equation}
From the functional limit theorem in Bojdcki et al (2006a), and assuming the one for equilibrium, a relationship analogous to (\ref{eq:4.1}) would hold for the two corresponding measure-valued limit processes (the constant $\kappa$ is the same in both cases). This raises the question of interpreting this mysterious result in terms of the particle system. Again going further, we think that for $\beta<1$ an analogous situation would appear. fBm would be replaced by the fractional stable process mentioned above, and sub-fBm by the sub-fractional stable process obtained in Bojdeki et al (2005), and they would be related analogously as in (\ref{eq:4.1}). (The sub-fractional stable process with $\beta =1$ is sub-fBm, Bojdecki et al, 2005; likewise, the fractional stable process with $\beta=1$ should be fBm). Then the  same question of interpretation would arise for the relationship between the equilibrium and Poisson fluctuation limit processes with $\beta <1$ in terms of the particle system.

Regarding the model of branching random walks on the lattice 
(Birkner and Z\"ahle, 2005), modifying it with the step of the walk in the domain of attraction of a spherically  symmetric $\alpha$-stable law, $(1+\beta)$-branching, $\beta <1$, and initial Poisson, the results are expected to be analogous to those in Bojdecki et al (2005) and the present paper.

\subsection{Long range dependence vs. independent increments}
\label{sub:4.4}

Clan recurrence and clan transience have been studied in 
by Stoeckel and Wakolbinger (1994) for $d\geq 3$, $\alpha =2$, $\beta =1$, in equilibrium. A clan is a family of infinitely many particles with eventually backwards coalescing paths, the system having started at time $-\infty$. It is shown that for $d=3$ (intermediate dimension) and $d=4$ (critical dimension), a.s. all clans visit every fixed ball in $\erre^d$ infinitely often and at arbitrarily large times (this is ``clan recurrence''), and for $d\geq 5$ (large dimensions) a.s. each clan visits a fixed ball in $\erre^d$ in a finite random interval of time and  never returns (this is ``clan transience''). Clan recurrence/transience for super-Brownian motion in equilibrium is discussed by Dawson and Perkins (1999). Although clan recurrence/transience has not been studied for the  branching particle system out of equilibrium (with $\alpha =2$, $\beta =1$),  one may think intuitively that the long range dependence of the occupation time fluctuation limit for $d=3$ is due to something close to clan recurrence: each one of many very large families of related particles visits  a fixed ball many times, each time adding   a random amount  to the occupation time of the ball. However, for $d=4$ this phenomenon still occurs but does not cause long range dependence. On the other hand, for $d\geq 5$ the independence of increments of the occupation time fluctuation limit (and the classical central limit norming $T^{1/2}$) may be attributed intuitively to something close to clan transience, since the  large families independently visit a fixed ball only up to  finite random  times. In order to explore these ideas it is necessary to formulate rigorously clan recurrence/transience out of equilibrium (at least with initial Poisson), also  for the general $(d, \alpha, \beta)$-branching particle system. 

In the case of intermediate dimensions with $\beta<1$ (Bojdecki et al 2005), an intriguing question is why there are two types of long range dependence regimes, one for $\beta>d/(d+\alpha)$ and another one for $\beta\leq d/(d+\alpha)$.

\subsection{Other long range dependence processes}
\label{sub:4.5}

For intermediate dimensions and $\beta =1$, Poisson initial condition leads to fBm and equilibrium initial condition leads to sub-fBm. The question here is if other initial conditions may lead to other long range dependence processes, different from fBm and sub-fBm. As a simple example, putting together two independent $(d, \alpha, 1)$-branching particle systems, one in equilibrium and the other one with initial Poisson (the normings are the same), the occupation time fluctuation limit would be the process $K\lambda (\xi + \zeta)$, where $\xi$  is fBm, 
$\zeta$ is sub-fBm, and they are independent (by Bojdecki et al, 2004, the constant $K$ is the same in both cases). Other, more interesting, examples will be considered elsewhere. The same kind of question arises for $\beta <1$. A large class of initial conditions (which includes Poisson), under which $N_t$ tends to equilibrium as $t\rightarrow \infty$, is considered by Gorostiza and Wakolbinger (1994).

Other long range dependence processes may be obtained by incorporating immigration in the branching particle systems. Some examples are given by Gorostiza et al (2005) without functional convergence proofs. In this case the fluctuations of the rescaled empirical process itself may have a long range dependence Gaussian limit process in dimension $d=1$ with $\alpha>1$ (see the covariance (1.5) in Gorostiza et al, 2005;  see also Li and Shiga, 1995, Theorem 1.10, for a superprocess setup). A surprising feature of this result is that the fluctuation of the rescaled empirical process is Markovian but its limit is not. The Markov property is not necessarily preserved under weak limits, but it would be interesting to understand why it is not preserved in this case.

\subsection{Superprocesses and other branching particle systems}
\label{sub:4.6}

Iscoe (1986) studied the occupation time  of the $(d, \alpha, \beta)$-superprocess, which is a measure-valued limit of the $(d, \alpha, \beta)$-branching system (see Dawson, 1993, Dawson and Perkins, 1991, Etheridge, 2000, for that class of superprocesses). He proved single time limit theorems for the occupation time fluctuations in large dimensions, $d>\alpha (1+\beta)/\beta$, and also for the intermediate dimension $d=3$ in the case $\alpha =2$, $\beta =1$. The limits are ${\cal S}'(\erre^d)$-valued random variables. For $\beta <1$ and large dimensions the limit is an ${\cal S}'(\erre^d)$-valued $(1+\beta)$-stable random variable, and a Hilbert subspace of ${\cal S}'(\erre^d)$ is found where this random field lives. (On this point it is relevant to mention that the aim of 
P\'erez-Abreu et al, 2005, is to determine, for a class of additive processes in the dual of a nuclear Fr\'echet space, a Hilbert subspace where such a process lives). For the intermediate dimension $d=3 \;(\alpha =2, \beta =1)$ an occupation time fluctuation limit process is obtained, with  fixed $\varphi \in {\cal S}(\erre^d)$, $\varphi \geq 0$, where the covariance function of sub-fBm appears with $h=3/2$; that was originally our motivation for investigating long range dependence or its absence in occupation time fluctuations of branching  systems.

Occupation time fluctuations for more general branching particle systems, with one and two levels of critical branching, have been studied by Dawson et al (2001). The results are single time limits. The atypical normings in the occupation time fluctuations for intermediate and critical dimensions with one branching level and $\beta=1$ arise as follows. Let $G_t \varphi = \int^t_0 {\cal T}_s \varphi ds$, $\varphi \geq 0$, $\varphi \not\equiv 0$ ($G_t \varphi \rightarrow \infty$ as $t\rightarrow \infty$ for transient particle motion). If $G^2_t \varphi$ grows like a function $f_t$ as $t\rightarrow \infty$, then the norming $F_T$ for the occupation time fluctuation limit has the form $(\int^T_0 f_s ds)^{1/2}$. For the standard $\alpha$-stable  process on $\erre^d$, $G^2_t$ grows like $t^{2-d/\alpha}$ for $\alpha <d<2\alpha$, and like log$t$ for $d=2\alpha$. Functional limit theorems for the branching systems in Dawson et al (2001) have not been attempted, except in the special cases considered in this paper and its 
predecessors.
\vglue .5cm
\noindent
{\large\bf Acknowlegments.} The authors thank the hospitality of the Institute of Mathematics, National University of Mexico (UNAM), where this paper was 
written. This research was partially supported by CONACyT grant 45684-F (Mexico) and MNiI grant 1P03A01129 (Poland).

\vglue1cm
\noindent
 Birkner M. and Z\"ahle, I. (2005) Functional limit theorems for the occupation time of the origin for branching random walks in $d\geq 3$, 
{\it Weierstrass Institut f\"ur Angewandte Analysis and Stochastik}, Berlin, preprint No. 1011.\\\\
 Bojdecki, T., Gorostiza, L.G. and Ramaswany, S. (1986) Convergence of ${\cal S}'$-valued processes and space-time random fields, {\it J. Funct. Anal}. 66, 21-41.\\\\
 Bojdecki, T., Gorostiza, L.G. and Talarczyk, A. (2004) Sub-fractional
Brownian motion and its relation to occupation times, {\it Stat. Prob.
Lett.} 69, 405-419.\\\\
 Bojdecki, T., Gorostiza, L.G. and Talarczyk, A. (2006a) Limit theorems for
occupation time fluctuations of branching systems I: long-range dependence,
{\it Stoch. Proc. Appl.} 116, 1-18.\\\\
 Bojdecki, T., Gorostiza, L.G. and Talarczyk, A. (2006b) Limit theorems for
occupation time fluctuations of branching systems II: critical and large
dimensions, {\it Stoch. Proc. Appl.} 116, 19-35.\\\\
 Bojdecki, T., Gorostiza, L.G. and Talarczyk, A. (2005) A long range dependence stable process and an infinite variance branching system, preprint.\\\\
 Dawson, D.A. (1993) Measure-Valued Markov Processes, Ecole d'\'et\'e de probabilit\'es de Saint-Flour XXI-1991, {\it Lect. Notes Math.} 1541, Springer.\\\\
 Dawson, D.A., Gorostiza, L.G., and Wakolbinger, A. (2001) Occupation time fluctuations in branching systems, {\it J. Theor. Probab.} 14, 729-796.\\\\
 Dawson, D.A. and Perkins, E. (1991) Historical Processes, {\it Memoirs Amer. Math. Soc.} 93, No. 454.\\\\
 Dawson, D.A. and Perkins, E. (1999) Measure-valued processes and renormalization of branching particle systems, in ``Stochastic Partial Differential Equations: Six Perspectives'', (Carmona, R. and Rozovskii, B., eds.), {\it Math. Surveys and Monographs} 64, 45-106, { Amer. Math. Soc.}\\\\
 Etheridge, A.M. (2000) An Introduction to Superprocesses, {\it University Lect. Series} 20, Amer. Math. Soc.\\\\
 Fleischmann, K. and G\"artner, J. (1986), Occupation time processes at a critical point, {\it Math. Nachr}. 125, 275-290.\\\\
 Gorostiza, L.G., Navarro, R. and Rodrigues, E.R. (2005), Some long-range dependence processes arising from fluctuations of particle systems, {\it Acta Appl. Math.} 86, 285-308.\\\\
 Gorostiza, L.G. and Wakolbinger, A. (1991) Persistence criteria for a class of critical branching particle systems in continuous time, {\it Ann. Probab}.  19, 266-288.\\\\
  Gorostiza, L.G. and Wakolbinger, A. (1994) Long time behaviour of critical branching particle systems and applications, in ``Measure-Valued Processes, Stochastic Partial Differential Equations, and Interacting Systems'', (Dawson, D.A., ed.), {\it CRM Proc. \&  Lect. Notes} 5,  119-137, Amer. Math. Soc.\\\\
 Iscoe, I. (1986) A weighted occupation time for a class of measure-valued branching processes, {\it Probab. Th. Rel. Fields} 71, 85-116.\\\\
 It\^o, K. (1980) Continuous additive ${\cal S}'$-processes, {\it Lect. Notes in Control and Inform. Sci.} 25, 36-46, Springer.\\\\
 It\^o, K. (1984) Foundations of Stochastic Differential Equations in Infinite Dimensional Spaces, {\it SIAM,} Philadelphia.\\\\
 Jakubowski, A. (1997) A non-Skorohod topology on the Skorohod space, {\it  Electron. J. Probab.} 2, no. 4, 21 pp.\\\\
 Kallenberg, O. (2002) Foundations of Moderm Probability, Second Edition,  Springer.\\\\
 Li, Z. and Shiga, T. (1995), Measure-valued diffusions: Immigrations, excursions and limit theorems, {\it J. Math. Kyoto Univ.} 35, 233-274.\\\\
 M\'el\'eard, S. and Roelly, S. (1992), An ergodic result for critical spatial branching processes, in ``Stochastic Analysis and Related Topics'' (Silivri, 1990), 333-341, {\it Progr. Probab.} 31, Birkh\"auser.\\\\
 P\'erez-Abreu, V., Rocha-Arteaga, A. and Tudor, C. (2005) Cone-additive processes in duals of nuclear Fr\'echet spaces, {\it Random Operators and Stochastic Evolution Equations} 13, 353-368.\\\\
 Skorohod, A.V. (1956) Limit theorems for stochastic processes, {\it Theory Probab. Appl.} 1, 261-290.\\\\
 Stoeckel, A. and Wakolbinger, A. (1994) On clan-recurrence and -transience in time stationary branching Brownian particle systems, in ``Measure-Valued Processes, Stochastic Partial Differential Equations, and Interacting Systems'', (Dawson, D.A., ed.), {\it CRM Proc. \& Lect. Notes} 5,  213-219, Amer. Math. Soc.\\\\
 Talarczyk, A. (2005) A functional ergodic theorem for the occupation time process of a branching system, preprint.\\\\
 \"Ustunel, A.S. (1984) Additive processes on nuclear spaces, {\it Ann. Probab.} 12, 858-868.
\end{document}